\documentclass[11pt]{amsart}
\usepackage{amsmath,amscd,amsthm,amssymb,amsxtra,latexsym,epsfig,epic,graphics}
\usepackage[matrix,arrow,curve]{xy}
\usepackage{graphicx}
\usepackage{MnSymbol}
\usepackage{tikz-cd}
\usepackage{color} 
\usetikzlibrary{arrows,calc} 
\usetikzlibrary{decorations.pathmorphing}
\usetikzlibrary{decorations.markings} 
\usetikzlibrary{decorations.pathreplacing} 
\usetikzlibrary{plothandlers}

\def\antiddot{\mathinner{\mkern1mu\raise1pt\vbox{\kern7pt\hbox{.}}\mkern2mu
        \raise4pt\hbox{.}\mkern2mu\raise7pt\hbox{.}\mkern1mu}}

\newcommand{\CC}{{\mathbb C}}

\newcommand{\NN}{{\mathbb N}}
\newcommand{\PP}{{\mathbb P}}
\newcommand{\QQ}{{\mathbb Q}}

\newcommand{\ZZ}{{\mathbb Z}}

\newcommand{\coker}{{\rm{coker}\,}}

\newcommand{\sF}{{\s F}}
\newcommand{\sG}{{\s G}}

\newcommand{\sJ}{{\s J}}

\newcommand{\sL}{{\s L}}

\newcommand{\sO}{{\s O}}

\newcommand{\punkt}{\hspace{-.3ex}\raise.15ex\hbox to1ex{\Huge.}}

\DeclareMathOperator{\Pic}{Pic}

\DeclareMathOperator{\Proj}{Proj}
\DeclareMathOperator{\Hom}{Hom}
\DeclareMathOperator{\sHom}{\sH om}

\DeclareMathOperator{\rank}{rk}

\newtheorem{theorem}{Theorem}[section]

\newtheorem{proposition}[theorem]{Proposition}

\theoremstyle{definition}

\def\PP{{\mathbb P}}

\def\QQ{{\mathbb Q}}

\def\CC{{\mathbb C}}

\def\sF{{\mathcal F}}
\def\sG{{\mathcal G}}
\def\sL{{\mathcal L}}
\def\sO{{\mathcal O}}
\def\sJ{{\mathcal J}}

\def\sHom{{\mathcal Hom}}

\def\sExt{{\mathcal Ext}}
\def\Hom{{\rm Hom}}

\DeclareMathOperator{\Eff}{{\rm Eff}}
\DeclareMathOperator{\Nef}{{\rm Nef}}
\DeclareMathOperator{\Mov}{{\rm Mov}}
\def\lbracket{{[\kern-1.5pt[}}
\def\rbracket{{]\kern-1.5pt]}}

\makeatletter
\def\Ddots{\mathinner{\mkern1mu\raise\p@
\vbox{\kern7\p@\hbox{.}}\mkern2mu
\raise4\p@\hbox{.}\mkern2mu\raise7\p@\hbox{.}\mkern1mu}}
\makeatother

\usepackage{times}
\newdimen\x \x=12pt

\usepackage{color}

\usepackage{hyperref}
	\author{Vladimir Lazi\'c}
	\address{Fachrichtung Mathematik, Campus, Geb\"aude E2.4, Universit\"at des Saarlandes, 66123 Saarbr\"ucken, Germany}
	\email{lazic@math.uni-sb.de}

\author{Frank-Olaf Schreyer}
\address{Fachrichtung Mathematik, Campus, Geb\"aude E2.4, Universit\"at des Saarlandes, 66123 Saarbr\"ucken, Germany}
\email{schreyer@math.uni-sb.de}

\title[Birational geometry of a family of determinantal 3-folds]{Birational geometry and the canonical ring\\ of a family of determinantal 3-folds}

\dedicatory{Dedicated to Giorgio Ottaviani on the occasion of his 60th birthday}

	\thanks{
Lazi\'c was supported by the DFG-Emmy-Noether-Nachwuchsgruppe ``Gute\linebreak Strukturen in der h\"oherdimensionalen birationalen Geometrie''. This work is also a contribution to Project-ID 286237555 -- TRR 195 -- by the Deutsche Forschungsgemeinschaft (DFG, German Research Foundation).
	}

\begin{document}

\begin{abstract} 
Few explicit families of 3-folds are known for which the computation of the canonical ring is accessible and the birational geometry non-trivial. In this  note we investigate a family of determinantal 3-folds in $\PP^2 \times \PP^3$ where this is the case.
\end{abstract}
	\maketitle
	\setcounter{tocdepth}{1}
	\tableofcontents

\section{Introduction}

There has been substantial progress in higher dimensional birational geometry over $\CC$ in the past decade. For instance, we currently know that for every smooth projective variety $X$, the canonical ring 
$$R(X,K_X)=\bigoplus_{m\in\NN}H^0(X,mK_X)$$ 
is finitely generated and that varieties with mild singularities and of log general type have good minimal models \cite{BCHM,CL12a,CL13}. Numerous other results have also recently been obtained when $X$ is not necessarily of general type, but the existence of minimal models and the Abundance conjecture remain unproven in general.

Lack of examples in higher dimensional geometry is one of the problems in the field for two reasons: (a) ultimately, one wants to apply the general theory in concrete examples, preferably described by concrete equations, and (b) without examples, it is often difficult to decide if a certain conjecture is plausible or to devise a route to a possible proof of a conjecture.

Recall that some of the main examples of higher dimensional constructions are the following: projective bundles (this is probably the most common class of examples, see \cite[2.3.B]{Laz04}); toric bundles, see \cite[Chapter IV]{Nak04}; deformations. Recently, blowups of $\mathbb P^3$ along a very general configuration of points were used in \cite{Les15} to give counterexample to  a conjecture of Kawamata, and a relatively simple example from \cite{Ogu14} (a complete intersection of general hypersurfaces of bi-degrees $(1,1)$, $(1,1)$ and $(2,2)$ in $\mathbb P^3\times\mathbb P^3$) was used in \cite{Les22} to disprove a widely believed claim from \cite{Nak04,Leh13,Eck16} about an expected behaviour of the numerical dimension.

The last two examples above should illustrate that more examples are needed in order to speed up progress in the field. We provide a general class of new examples in this note, and investigate the birational geometry of a particular subclass of examples in detail.

\medskip

The class of examples we study in this paper are a particular case of \emph{determinantal varieties}. The situation in general is explained in detail in Section \ref{sec:2}. In particular, denote $\PP=\PP^2\times \PP^3$ and $\sF = \sO_\PP^{\oplus2}$, and for each integer $b \ge 1$ define the sheaf
$$\sG_b = \sO_\PP(1,b)\oplus \ker\big(H^0(\PP,\sO_\PP(1,0)) \otimes \sO_\PP(1,1) \to \sO_\PP(2,1)\big).$$
Pick $\varphi \in \Hom(\sF,\sG_b)$ general, and $X_b$ let be the $3$-fold given as
$$X_b=\{p\in \PP \mid \rank \varphi(p) \leq 1\}.$$
Our main result is:

\begin{theorem}\label{thm:main} 
The variety $X_b$ is birational to a hypersurface $Y_b$ of degree $ 2b+2$ in the weighted projective space $\PP(1,1,1,1,b+1)$. In particular, we have
$$\kappa(X_b) =
\begin{cases}
-\infty &\hbox{ if }\; b=1\hbox{ or } 2,\cr
0 &\hbox{ if }\; b=3,\cr
3  &\hbox{ if }\; b\ge 4.
\end{cases}
 $$ 
The image $X_b^1$ of $X$ in $\PP^1\times \PP^3$ is a small resolution of $Y_b$ in  $(b+1)^3$ $A_1$-sin\-gu\-la\-ri\-ties.  The morphism $X \to X_b^1$ is the blowup of one of the two components of the preimage of a twisted $C \subseteq \PP^3$ which intersects the branch divisor of $Y_b \to \PP^3$ tangentially. The variety $X_b^1$ has precisely two minimal models and one nontrivial birational automorphism $\iota$ of order two. The automorphism $\iota$ interchanges the two models.
\end{theorem}

Thus for $b\geq4$ the 3-fold $X_b^1$ is a minimal model of $X$ and $Y_b$ is the canonical model. In particular, this family of examples has an unexpectedly rich birational geometry.

\section{Determinantal varieties}\label{sec:2}

In this section we describe a general construction of determinantal varieties in products of projective spaces, and specialise to a particular case which is the main object of this paper.

\subsection{A general construction}

Let $\PP$ be a product of projective spaces, let $\sF$ and $\sG$ be vector bundles on $\PP$ of rank $f$ and $g \ge f$ respectively, and let $\varphi\colon \sF \to \sG$ be a general homomorphism. Define an algebraic set $X\subseteq\PP$ by
$$X = \{ p \in \PP \mid \varphi(p) \hbox{ does not have maximal rank } f \}.$$

For example, if the sheaf $\sHom_{\sO_\PP}(\sF,\sG)\simeq\sF^*\otimes_{\sO_\PP}\sG$ is ample, then $X$ is non-empty, connected and has codimension $g-f+1$ by \cite{FL81}, and is smooth outside a sublocus of codimension $2(g-f+2)$ by \cite{Kle69}, which is empty if $\dim \PP < 2(g-f+2)$. Moreover, in this case the sheaf
\begin{equation}\label{eq:1}
\sL= \coker(\varphi^t\colon  \sG^* \to \sF^*)
\end{equation}
is a line bundle on $X$. 

If $f=1$, then $X$ is a zero loci of a section of a vector bundle on $\PP$. If additionally $\sG$ is a direct sum of line bundles, then $X$ is a complete intersection.

\medskip

Perhaps the simplest case beyond the one above is when $f=g-1$. In that case, $X$ is a codimension $2$ subvariety in $\PP$ and, if $\sJ_X$ is the ideal sheaf of $X$ in $\PP$, then we have the resolution
\begin{equation}\label{eq:2}
0 \to \sF \, {\buildrel {\varphi} \over {\longrightarrow} }\,  \sG \to \sJ_X\otimes\sO_\PP\big(c_1(\sG)-c_1(\sF)\big) \to 0,
\end{equation}
see \cite{DES93,DP95}. By above, we expect these $X$ to be a smooth variety only when $\dim \PP \le 5$.

\subsection{Examples}

Thus, from now on we choose $\PP=\PP^2 \times \PP^3$, and we let $f=2$ and $g=3$. Specifying further $\sF:= \sO_\PP^{\oplus2}$, then $X$ is a $3$-fold and the linear system $|\sL|$, where $\sL$ is defined as in \eqref{eq:1}, defines a morphism $\PP\to\PP^1$. Since we also have the projections from $\PP$ to its two factors, we obtain three maps
\begin{equation}\label{eq:2b}
\pi_1\colon  X \to \PP^1, \quad \pi_2\colon  X \to \PP^2, \quad \pi_3\colon  X \to \PP^3,
\end{equation} 
which we use to study $X$.

\medskip

At first sight, the case $\sG=\sO_\PP(1,1)^{\oplus3}$ might look like the simplest possible case. In this case, the morphism $\pi_2\colon  X \to \PP^2$ is a fibration into twisted cubic curves, $\pi_3\colon  X \to \PP^3$ is generically finite of degree $3:1$, and $\pi_1\colon X \to \PP^1$ is a fibration into cubic surfaces. 

\medskip

Now, let $\theta\colon \sO_\PP(1,1)^{\oplus4} \to \sO_\PP(2,1)$ be a general morphism and consider the case $\sG= \ker\theta$. In suitable coordinates on $\PP^2$ we have
$$\mathcal G=\sO_\PP(1,1)\oplus \ker\big(H^0(\PP,\sO_\PP(1,0)) \otimes \sO_\PP(1,1) \to \sO_\PP(2,1)\big),$$
where the map is the evaluation morphism. This case is even simpler, in the sense that $\pi_3\colon  X \to \PP^3$ is generically finite of degree $2:1$. Indeed, let $F$ be a general fiber of the second projection $\PP \to \PP^3$. Then the sheaf
$$ \sG|_F\simeq \ker\big(\sO_{\PP^2}(1)^{\oplus4} \to \sO_{\PP^2}(2)\big)$$
has the Chern polynomial
$$ c_t(\sG|_F)= \frac{(1+t)^4}{1+2t}=1+2t+2t^2,$$
and thus $c_2(\sG|_F)=2$ implies that $\pi_3$ is generically $2:1$.

\section{Cohomological properties}

\subsection{The main example}
Our main example is a generalisation of this last construction. As announced in the introduction, for each integer $b \geq 1$ we consider $3$-folds $X_b$ constructed as follows: we set $\PP=\PP^2 \times \PP^3$, $\sF= \sO_\PP^{\oplus2}$, and
$$\sG_b = \sO_\PP(1,b)\oplus \ker\big(H^0(\PP,\sO_\PP(1,0)) \otimes \sO_\PP(1,1) \to \sO_\PP(2,1)\big),$$
where the morphism is the evaluation morphism in suitable coordinates $(x_0:x_1:x_2)$ on $\PP^2$. Then for a general $\varphi \in \Hom(\sF,\sG_g)$ we define
$$X_b=\{p\in \PP^2\times \PP^3 \mid \rank \varphi(p) \leq1\}.$$
This is the main object of this paper.

By \eqref{eq:2}, there exists a locally free resolution
\begin{equation}\label{ideal}
0 \to \sO_\PP^{\oplus2} \to \sG_b \to \sJ_{X_b}\otimes\sO_\PP(2,b+2) \to 0,
\end{equation}
and $ \pi_3 \colon X_b \to \PP^3$ is generically $2:1$ similarly as in Section \ref{sec:2}. Dualizing \eqref{ideal} we obtain a resolution of $\sL$:
\begin{align}\label{ell}
0 \leftarrow \sL \leftarrow \sO_\PP^{\oplus2} &\leftarrow  \sO_\PP(-1,-1)^{\oplus3} \oplus \sO_\PP(-1,-b)\\
& \leftarrow \sO_\PP(-2,-1) \oplus \sO_\PP(-2,-2-b) \leftarrow 0, \notag
\end{align}
and thus 
\begin{align}
\omega_{X_b} &\simeq \sExt^2_{\sO_\PP}\big(\sO_{X_b},\sO_\PP(-3,-4)\big) \label{eq:omega}\\
&\simeq \sExt^2_{\sO_\PP}\big(\sO_{X_b}(2,b+2),\sO_\PP(-1,b-2)\big) \simeq \sL(-1,b-2).\notag
\end{align}

Some of the computationally accessible information in explicit examples are the dimensions of the cohomology groups 
$H^i\big(X_b,\sO_{X_b}(\alpha,\beta)\big)$. It is useful to arrange this data in cohomology polynomials
$$p_{\alpha,\beta}= \sum_{i=0}^3 h^i\big(X_b,\sO_{X_b}(\alpha,\beta)\big) \cdot h^i \in \ZZ[h].$$ 
We also consider the ring
$$R=\bigoplus_{\beta \ge 0} H^0\big(X_b,\sO_{X_b}(0,\beta)\big).$$

\subsection{Cohomology groups of $X_3$}\label{subsec:X3}
Using the theory of Tate resolutions for product of projective spaces \cite{EES15} we can calculate the dimensions of these groups. In this subsection, we concentrate on the case $b=3$. Fix the range
$$ -3 \le \alpha \le 3, \quad -7 \le \beta \le 7.$$
Then we can summarize the result in matrix of cohomology polynomials $p_{\alpha,\beta}$ as below.

\noindent
{\small \begin{tabular}{rrrrrrr}
      $88 h$& $
      56 h$& $
      20$& $
      140$& $
      304$& $
      512$& $
      764$\\$
      53 h$& $
      41 h$& $
      8$& $
      94$& $
      217$& $
      377$& $
      574$\\$
      24 h$& $
      26 h$& $
      2$& $
      60$& $
      148$& $
      266$& $
      414$\\$
      5 h^{2}+8 h$& $
      13 h$& $
      0$& $
      36$& $
      95$& $
      177$& $
      282$\\$
      10 h^{2}+2 h$& $
      4 h$& $
      0$& $
      20$& $
      56$& $
      108$& $
      176$\\$
      7 h^{2}$& $
      h$& $
      0$& $
      10$& $
      29$& $
      57$& $
      94$\\$
      12 h^{3}+4 h^{2}$& $
      6 h^{3}$& $
      2 h^{3}$& $
      4$& $
      12$& $
      2 h+24$& $
      6 h+40$\\$
      40 h^{3}+h^{2}$& $
      21 h^{3}$& $
      8 h^{3}$& $
      h^{3}+1$& $
      3$& $
      5 h+6$& $
      16 h+10$\\$
      88 h^{3}$& $
      48 h^{3}$& $
      20 h^{3}$& $
      4 h^{3}$& $
      0$& $
      8 h$& $
      28 h$\\$
      157 h^{3}$& $
      89 h^{3}$& $
      40 h^{3}$& $
      10 h^{3}$& $
      h^{2}$& $
      h^{2}+8 h$& $
      34 h$\\$
      248 h^{3}$& $
      146 h^{3}$& $
      70 h^{3}$& $
      20 h^{3}$& $
      4 h^{2}$& $
      4 h^{2}+2 h$& $
      2 h^{2}+28 h$\\$
      363 h^{3}$& $
      221 h^{3}$& $
      112 h^{3}$& $
      36 h^{3}$& $
      7 h^{2}$& $
      17 h^{2}$& $
      8 h^{2}+14 h$\\$
      504 h^{3}$& $
      316 h^{3}$& $
      168 h^{3}$& $
      60 h^{3}$& $
      2 h^{3}+10 h^{2}$& $
      36 h^{2}$& $
      24 h^{2}$\\$
      673 h^{3}$& $
      433 h^{3}$& $
      240 h^{3}$& $
      94 h^{3}$& $
      8 h^{3}+13 h^{2}$& $
      57 h^{2}$& $
      62 h^{2}$\\$
      872 h^{3}$& $
      574 h^{3}$& $
      330 h^{3}$& $
      140 h^{3}$& $
      20 h^{3}+16 h^{2}$& $
      78 h^{2}$& $
      106 h^{2}$\\
      \end{tabular}}

Let us point out a few interesting values: we have
$$h^1(X_b,\sO_{X_b})=h^2({X_b},\sO_{X_b})=0 \hbox{ and }h^3({X_b}, \sO_{X_b})=h^0({X_b}, \omega_{X_b})=1 $$
from the center entry. Moreover, we see that $h^0\big(X_b,\sO_{X_b}(0,4)\big)=36 > 35$, so the ring $R$ has a further generator in degree $4$.

Another interesting sequence of values are the dimensions of the $H^2$-co\-ho\-mo\-lo\-gy in the first vertical strand (that is, for $\alpha=1$): 
$$\ldots,16,13,10,7,4,1.$$
This looks like the Hilbert function of  the twisted cubic in $\PP^3$.

\subsection{Cohomology groups of $X_b$}
The tables for other values of $b$ have a lot of similarity with the table above.

Recall from \eqref{eq:omega} that $\sL\cong \omega_{X_b}(1,-b+2)$. Dualising the resolution \eqref{ell} we obtain
\begin{align*}
0 &\to \sO_\PP^{\oplus2} \to \sO_\PP(1,1)^{\oplus3}\oplus \sO_\PP(1,b)\\
& \to \sO_\PP(2,1) \oplus \sO_\PP(2,b+2) \to \sO_{X_b}(2,b+2) \to 0.
\end{align*}
Twisting back by $\sO_\PP(-2,-b-2)$ we deduce
$$R\pi_{3,*} \sO_{X_b} = \pi_{3,*} \sO_{X_b} = \sO_{\PP^3}\oplus \sO_{\PP^3}(-b-1),$$
and twisting by back by $\sO_\PP(-3,-b-2)$ gives
$$R\pi_{3,*} \sO_{X_b}(-1,0) = \sO_{\PP^3}(-b-2)^{\oplus2}.$$
Since $R\pi_{3,*}\sO_{X_b}(\alpha,0)$ is computed with the vertical strands in the Tate resolution, this explains the values in the $0$-th and $(-1)$-st vertical strand in the cohomology table. In particular, we see that 
$$h^0\big(X_b,\sO_{X_b}(-1,b+2)\big)=2.$$

\subsection{A twisted cubic}\label{subsec:twisted}
As suggested in \S\ref{subsec:X3}, we can find a twisted cubic on $\PP^3$ in our construction.

Recall that we fixed coordinates $(x_0:x_1:x_2)$ on $\PP^2$. We may write
	\begin{center}
		\begin{tikzcd}
			\sG_b = \sO_\PP(1,b)\oplus \ker\big(\sO_\PP(1,1)^{\otimes3} \arrow[rr, "{(x_0,x_1,x_2)}"] && \sO_\PP(2,1)\big)
		\end{tikzcd}
	\end{center} 

so that we have two projections
$$\sG_b\to\sO_\PP(1,b)\quad\text{and}\quad \sG_b\to\sO_\PP(1,1)^{\otimes3}.$$
The composition
$$ \sO_\PP^{\oplus2}  \, {\buildrel {\varphi} \over {\longrightarrow} }\,  \sG_b \to \sO_\PP(1,1)^{\oplus3}$$
factors over
$$\sO_\PP^{\oplus2} \, {\buildrel {\psi} \over {\longrightarrow} }\,  \sO_\PP(0,1)^{\oplus3} \, {\buildrel {K_2} \over {\longrightarrow} }\,  \sO_\PP(1,1)^{\oplus3},$$
where
$$
K_2=\begin{pmatrix}
 0 & -x_2 & x_1 \cr
 x_2 & 0 & -x_0 \cr
 -x_1 & x_0 & 0 \cr
\end{pmatrix}
$$
is the Koszul matrix, and in suitable coordinates $(y_0:y_1:y_2:y_3)$ of $\PP^3$ we have
$$
\psi = \begin{pmatrix}
 y_0 & y_1 & y_2 \cr
 y_1 & y_2& y_3 \cr
\end{pmatrix}.
$$
We denote by $C \subseteq \PP^3$ the twisted cubic curve defined by the $2\times 2 $ minors of $\psi$.

\medskip

The remaining part $\sO_\PP^{\oplus2} \to \sO_\PP(1,b)$ of $\varphi$ can be factored as $B \cdot \begin{pmatrix}
x_0 & x_1 &x_2
\end{pmatrix}^t$, with
\begin{equation}\label{eq:B}
B=\begin{pmatrix} 
b_{00} & b_{01} & b_{02} \cr
b_{10} & b_{11} & b_{12} \cr
\end{pmatrix},
\end{equation}
where $b_{ij}\in \CC[y_0,y_1,y_2,y_3]$ are forms of degree $b$. To this matrix we associate the matrix
\begin{equation}\label{eq:M}
M=\begin{pmatrix}2 
\sum\limits_{i=0}^2 y_ib_{0i} & \sum\limits_{i=0}^2 y_ib_{1i}+\sum\limits_{i=0}^2 y_{i+1}b_{0i} \cr
 \sum\limits_{i=0}^2 y_ib_{1i}+\sum_{i=0}^2 y_{i+1}b_{0i}   & 2\sum\limits_{i=0}^2 y_{i+1}b_{1i}  
\end{pmatrix};
\end{equation}
this matrix will be important in \S\ref{subsec:Xb1} below.

\begin{proposition}\label{pro:twisted}
In the notation as above, we have:
\begin{enumerate}
\item[(a)] $\pi_3^{-1}(C)\subseteq X_b$ decomposes into two components: $C_1$ of dimension $1$ and $E$ of dimension $2$,
\item[(b)] $C_1$ is defined by the $2 \times 2$ minors of
$$
\begin{pmatrix}
y_0 & y_1 & y_2 & x_0 & x_1\cr
 y_1 & y_2 & y_3 & x_1 & x_2 \cr
\end{pmatrix},
$$
\item[(c)] $E$ is defined by the minors of $\psi$ and the entries of
$$
\begin{pmatrix}
x_0 & x_1 &x_2
\end{pmatrix}
\cdot B^t \cdot \begin{pmatrix} 
0 & 1 \cr 
-1 & 0 \cr
\end{pmatrix}
\cdot \psi,
$$
\item[(d)] $C_1 \to C$ is an isomorphism while $E \to C$ is a $\PP^1$-bundle. In particular, $C_1$ and $E$ are smooth.
\end{enumerate}
\end{proposition}

\begin{proof}
Parts (b) and (c) follow from direct calculations \cite{MacPac} or \cite{LS19}. Note that
$$\{ p \in \PP^3 \mid \rank B(p)\le 1 \hbox{ and } \rank \psi(p) \le 1 \} = \emptyset$$ 
for a general choice of $B$. Therefore, $\rank B(p)=2$ for $p\in C$, so $B^t \cdot \begin{pmatrix} 
0 & 1 \cr 
-1 & 1 \cr
\end{pmatrix}
\cdot \psi$
has rank $1$ over the points of $C$. Hence, $E$ is a $\PP^1$-bundle. We have $C_1 \cong C$ and the projection $\pi_2$ maps $C_1$ isomorphically to the conic $V(x_0x_2-x_1^2) \subseteq \PP^2$. This shows (d).

Finally, consider the matrix $\varphi^t$ as a $2 \times 4$ matrix with entries in 
$$\QQ[x_0,x_1,x_2,y_0,y_1,y_2,y_3,b_{00},\ldots,b_{12}].$$ 
The defining ideal of $X_b$ is the annihilator of the $\coker \varphi$, once we substitute the actual values for the $b_{ij}$  in $H^0\big(\PP,\sO_\PP(0,b)\big)$. Adding the defining equations of $C$, a primary decomposition gives the two components in this generic setting. Since $C_1$ and $E$ are smooth, specialising $b_{ij}$ gives the actual components. \end{proof}

\section{Two minimal models}

In this section we describe the birational geometry of $X_b$. 

\subsection{An overview}\label{subsec:overview}

We introduce several new varieties. Denote
$$X_b^1:=(\pi_1 \times \pi_3)(X_b) \subseteq \PP^1 \times \PP^3.$$
Moreover, let 
$$
R=\CC[y_0,y_1,y_2,y_3,w]/\langle w^2+\det M \rangle,
$$
where $w$ has degree $b+1$ and $M$ is defined as in \eqref{eq:M}, and denote 
$$Y_b =  \Proj R  \subseteq \PP(1,1,1,1,b+1) .$$
An easy argument with an exact sequence in \S\ref{subsec:Xb1} shows the existence of a rational map $\rho\colon  X_b^1 \dasharrow \PP^1$, and we denote
$$X_b^2:=(\rho\times \pi_3)(X_b^1)\subseteq\PP^1 \times \PP^3.$$

We will show that these varieties fit into the diagram
\begin{equation}
	\begin{gathered}
		\begin{tikzcd}
			X_b \arrow[d, "\alpha_1" swap] \arrow[r, dashed, "\alpha_2"] & X_b^2 \arrow[d, "\xi_2"]  \\
			X_b^1 \arrow[ru, dashed] \arrow[r, "\xi_1" swap] & Y_b, 
		\end{tikzcd}
	\end{gathered} 
\end{equation}
such that the following holds:
\begin{enumerate}
\item[(a)] $X_b^1 \subseteq \PP^1\times \PP^3$ is a hypersurface of bi-degree $(2,b+1)$,
\item[(b)] $X_b^1$ and $X_b^2$ are small resolutions of $Y_b$.
\end{enumerate}
This then implies our main result.

\subsection{The geometry of $X_b^1$}\label{subsec:Xb1}

Our first goal is to compute $X_b^1$. 

By \S\ref{subsec:twisted}, the defining ideal of $X_b \subseteq \PP^1 \times \PP^2\times \PP^3$ is given by the four entries of the matrix
\begin{equation}\label{eq:Xb}
\begin{pmatrix} z_0 & z_1
\end{pmatrix}
\cdot
\big[\psi \cdot K_2\, \big| \, B \cdot \begin{pmatrix} x_0 & x_1 & x_2 \end{pmatrix}^t \big].
\end{equation}
The saturation of this ideal with respect to $\langle x_0,x_1,x_2 \rangle$ gives the hypersurface $X_b^1$. 

\begin{proposition}\label{pro:Xb1}
With notation as in \S\ref{subsec:overview}, we have:
\begin{enumerate}
\item[(a)] The variety $X_b^1$ is a smooth hypersurface of bi-degree $(2,b+1)$ in $\PP^1\times \PP^3$ defined by
$$
f=\begin{pmatrix} z_0 & z_1
\end{pmatrix}
\cdot M \cdot
\begin{pmatrix} z_0 \cr z_1
\end{pmatrix},
$$
with matrix $M$ given as in \eqref{eq:M}.
\item[(b)] The map $\alpha_1\colon X_b \to X_b^1$ is birational: it is the blow down of the $\PP^1$-bundle $E$ from Proposition \ref{pro:twisted} to the rational curve $C^1 \subseteq X_b^1$ defined by the $2 \times 2$ minors of the matrix
$$
\begin{pmatrix}
y_0 & y_1 & y_2 & -z_1 \cr
 y_1 & y_2 & y_3 & z_0 \cr
\end{pmatrix}.
$$
\end{enumerate}
\end{proposition}

\begin{proof}
We rewrite the equation \eqref{eq:Xb} of $X_b$ as
$$ 
\begin{pmatrix} x_0 & x_1 &x_2 
\end{pmatrix} \cdot N=0,$$
where
$$
N=\begin{pmatrix}0&
      {z}_{0} {y}_{2}+{z}_{1} {y}_{3}&
      -{z}_{0} {y}_{1}-{z}_{1} {y}_{2}&
      {z}_{0} {b}_{00}+{z}_{1} {b}_{10}\\
      -{z}_{0} {y}_{2}-{z}_{1} {y}_{3}&
      0&
      {z}_{0} {y}_{0}+{z}_{1} {y}_{1}&
      {z}_{0} {b}_{01}+{z}_{1} {b}_{11}\\
      {z}_{0} {y}_{1}+{z}_{1} {y}_{2}&
      -{z}_{0} {y}_{0}-{z}_{1} {y}_{1}&
      0&
      {z}_{0} {b}_{02}+{z}_{1} {b}_{12}\\
      \end{pmatrix}.
$$
We conclude that $X_b^1 \subseteq \PP^1 \times \PP^3$ coincides  with the variety defined by the radical of the $3\times 3$ minors of $N$. This radical coincides with the form $f$ in the statement of the proposition; the details of the calculations are in \cite{MacPac} or \cite{LS19}. Moreover, the map $\alpha_1$ is birational outside the preimage of the ideal defined by the $2 \times 2$ minors of $N$: this is the curve $C^1$. Since $\alpha_1$ blows down a smooth $\PP^1$-bundle $E$, the variety $X_b^1$ is smooth.
\end{proof}

With this information, one can calculate the cohomology table of $X_b^1$ in the test case $b=3$, using the Macaulay2 package TateOnProducts: 

\medskip
\noindent
{\small
\begin{tabular}{rrrrrrrrr}
    $148 h$& $
       96 h$& $
       44 h$& $
       8$& $
       60$& $
       112$& $
       164$& $
       216$& $
       268$\\ $
       100 h$& $
       66 h$& $
       32 h$& $
       2$& $
       36$& $
       70$& $
       104$& $
       138$& $
       172$\\ $
       60 h$& $
       40 h$& $
       20 h$& $
       0$& $
       20$& $
       40$& $
       60$& $
       80$& $
       100$\\ $
       30 h$& $
       20 h$& $
       10 h$& $
       0$& $
       10$& $
       20$& $
       30$& $
       40$& $
       50$\\ $
       12 h$& $
       8 h$& $
       4 h$& $
       0$& $
       4$& $
       8$& $
       12$& $
       16$& $
       20$\\ $
       5 h^{3}+3 h$& $
       4 h^{3}+2 h$& $
       3 h^{3}+h$& $
       2 h^{3}$& $
       h^{3}+1$& $
       2$& $
       h^{2}+3$& $
       2 h^{2}+4$& $
       3 h^{2}+5$\\ $
       20 h^{3}$& $
       16 h^{3}$& $
       12 h^{3}$& $
       8 h^{3}$& $
       4 h^{3}$& $
       0$& $
       4 h^{2}$& $
       8 h^{2}$& $
       12 h^{2}$\\ $
       50 h^{3}$& $
       40 h^{3}$& $
       30 h^{3}$& $
       20 h^{3}$& $
       10 h^{3}$& $
       0$& $
       10 h^{2}$& $
       20 h^{2}$& $
       30 h^{2}$\\ $
       100 h^{3}$& $
       80 h^{3}$& $
       60 h^{3}$& $
       40 h^{3}$& $
       20 h^{3}$& $
       0$& $
       20 h^{2}$& $
       40 h^{2}$& $
       60 h^{2}$\\ $
       172 h^{3}$& $
       138 h^{3}$& $
       104 h^{3}$& $
       70 h^{3}$& $
       36 h^{3}$& $
       2 h^{3}$& $
       32 h^{2}$& $
       66 h^{2}$& $
       100 h^{2}$\\ $
       268 h^{3}$& $
       216 h^{3}$& $
       164 h^{3}$& $
       112 h^{3}$& $
       60 h^{3}$& $
       8 h^{3}$& $
       44 h^{2}$& $
       96 h^{2}$& $
       148 h^{2}$\\
       \end{tabular}}

    \medskip
       
From the table, we see that $h^0\big(X_b^1,\sO_{X_b^1}(-1,4)\big)=2$. In fact, for every $b$ we have 
\begin{equation}\label{eq:b+1}
h^0 \big(X_b^1,\sO_{X_b^1}(-1,b+1)=2.
\end{equation}
This follows from the exact sequence
$$ 0 \to \sO_{\PP^1\times \PP^3}(-3,0) \to \sO_{\PP^1\times \PP^3}(-1,b+1)\to  \sO_{X_b^1}(-1,b+1) \to 0$$
and the fact that $h^1\big(\PP^1 \times \PP^3,\sO_{\PP^1 \times \PP^3}(-3,0)\big)=2$. Therefore, as announced in \S\ref{subsec:overview}, by \eqref{eq:b+1} we obtain a rational map
\begin{equation}\label{eq:rho}
\rho\colon  X_b^1 \dasharrow \PP^1.
\end{equation}

\subsection{The first small resolution}

Next we show that $X_b^1$ is a small resolution of $Y_b$ and analyse in detail the geometry of $Y_b$. Recall that by the definition of $Y_b$ in \S\ref{subsec:overview}, there exists a double cover 
\begin{equation}\label{eq:delta}
\delta\colon Y_b\to \PP^3.
\end{equation}

\begin{proposition}\label{pro:A1}
For a general choice of $b_{ij}$ in \eqref{eq:B} we have:
\begin{enumerate}
\item[(a)] the double cover $\delta$ has $A_1$-singularities above the  $(b+1)^3$ distinct points defined by the zero loci of entries of $M$, and is otherwise smooth,
\item[(b)] $X_b^1$ is a small resolution of $Y_b$.
\end{enumerate}
\end{proposition}

\begin{proof} 
Recall that the variety $X_b^1$ comes with a projection to $\PP^3$. By the description in Proposition \ref{pro:Xb1}, the fibre of the map $X_b^1 \to \PP^3$ over a point $p \in \PP^3$ consist either of two points, of one point or is isomorphic to $\PP^1$, depending on whether $M(p)$ has rank $2$, $1$ or $0$ respectively. For general $b_{ij}$, the three entries of the matrix $M$ form a regular sequence, which intersect in $(b+1)^3$ distinct points. Since this is an open condition for the values of $b_{ij}$, it suffices to construct an example. 

To this end, pick $\lambda_0, \ldots, \lambda_b, \mu_0, \ldots \mu_b \in \CC$ which are algebraically independent over $\QQ$. Define forms 
$$\widetilde b_{01}\in\QQ[\lambda_0, \ldots, \lambda_b][y_0,y_1],\qquad \widetilde{b}_{11}\in\QQ[\mu_0, \ldots \mu_b ][y_2,y_3]$$
of degree $b$ by the relations
$$ 
\prod_{i=0}^b (y_0-\lambda_i y_1)= y_0^{b+1}+y_1 \widetilde b_{01},\quad \prod_{j=0}^b (y_3-\mu_j y_2)= y_3^{b+1}+y_2 \widetilde b_{11},
$$
and define the matrix
$$
B^\circ= \begin{pmatrix} y_0^b &\widetilde b_{01} &0 \cr 0 &\widetilde b_{11} & y_3^b \cr \end{pmatrix}.
$$
We consider $B^\circ$ as the matrix $B$ from \eqref{eq:B} for special values of $b_{ij}$. For these values, the corresponding matrix $M$ from \eqref{eq:M} turns into
$$
M^\circ=\begin{pmatrix}
2\big(y_0^{b+1}+y_1 \widetilde b_{01}\big) & y_0^by_1+y_1\widetilde{b}_{11}+y_2y_3^b+y_2\widetilde{b}_{01} \cr
 y_0^by_1+y_1\widetilde{b}_{11}+y_2y_3^b+y_2\widetilde{b}_{01}   & 2\big(y_3^{b+1}+y_2 \widetilde b_{11}\big)  
\end{pmatrix}.
$$
Fix $0\leq i,j\leq b$. The diagonal entries of $M^\circ $ have solutions $y_0=\lambda_iy_1$ and $y_3=\mu_jy_2$. Substituting these values for $y_0$ and $y_3$ into the off diagonal entry of $M^\circ$ yields non-zero polynomials
\begin{align*}
P_{ij}&=\lambda_i^by_1^{b+1}+y_1\widetilde{b}_{11}(y_2,\mu_j y_2)+\mu_j^by_2^{b+1}+y_2\widetilde{b}_{01}(\lambda_iy_1,y_1) \\
&=\lambda_i^by_1^{b+1}-(\mu_j^{b+1}+ \ldots)y_1y_2^b+\mu_j^by_2^{b+1}-(\lambda_i^{b+1}+ \ldots)y_2y_1^b \\
&\in \QQ[\lambda_0,\ldots,\lambda_b,\mu_1,\ldots, \mu_b][y_1,y_2].
\end{align*}
The highest exponent of $\lambda_i$ and $\mu_j$  in the Sylvester matrix for the resultant
$$
\textstyle R\big(\frac{\partial P_{ij}}{\partial y_1} , \frac{\partial P_{ij}}{\partial y_2}\big)
$$ 
is $b+1$ and the coefficient of $(\lambda_i\mu_j)^{b(b+1)}$ is $\pm 1$ obtained from the coefficient of $y_2^b$ in $\frac{\partial P_{ij}}{\partial y_1}$ and the coefficient of $y_1^b$ in $\frac{\partial P_{ij}}{\partial y_2}$. Hence, the discriminant of $P_{ij}$ in $\QQ[\lambda_0,\ldots,\lambda_b,\mu_1,\ldots, \mu_b]$ is not identically zero. Since $\lambda_0,\ldots,\lambda_b,\mu_1,\ldots, \mu_b$ are algebraically independent over $\QQ$,  each $P_{ij}$ factors into $b+1$ distinct linear forms in $\CC[y_1,y_2]$. Hence, the entries of $M^\circ$ vanish in precisely $(b+1)^3$ distinct points, as desired.

Now, write 
$$M=\begin{pmatrix}
a_0 & a_1 \cr
a_1 & a_2 \cr
\end{pmatrix}
$$ 
for forms $a_i$ of degree $b+1$ on $\PP^3$ as in \eqref{eq:M}. For any $B$ leading to $(b+1)^3$ distinct points in $\PP^3$, the entries $a_0,a_1,a_2$ generate locally at each point its maximal ideal, so the branch divisor $\det M=0$ 
has $A_1$-singularities at these points. Since $X_b^1$ is smooth by Proposition \ref{pro:Xb1}, the branch divisor $\det M = 0$ is smooth outside the $A_1$-singularities.

Consider the subvariety of $\PP^1\times \PP(1,1,1,1,b+1)$ defined by the $2 \times 2$ minors of the matrix
\begin{equation}\label{eq:firstsmall}
\begin{pmatrix}
{a}_{0}&
{a}_{1}-w&
      {z}_{1}\\
{a}_{1}+w&
      {a}_{2}&      
      -{z}_{0}
      \end{pmatrix}.
\end{equation}
This is a small resolution of $Y_b$, and it is easy to see that it is isomorphic to $X_b^1$, as defined in Proposition \ref{pro:Xb1}(a).
\end{proof}

\begin{proposition}
Let $C\subseteq\PP^3$ be the twisted cubic defined in \S\ref{subsec:twisted} and let $\delta$ be the double cover from \eqref{eq:delta}. Then $C$ intersects the branch divisor of $\delta$ tangentially. We have $(\delta\circ\xi_1)^{-1}(C)=C^1\cup C^2\subseteq X_b^1$, where $C^1$ is the curve from Proposition \ref{pro:Xb1}, and $C^2$ is defined by the $4\times 4$ Pfaffians of the matrix
  $$\begin{pmatrix}0&
       0&
       {y}_{1}&
       {y}_{2}&
       {y}_{3}\\
       0&
       0&
       {y}_{0}&
       {y}_{1}&
       {y}_{2}\\
       {-{y}_{0}}&
       {-{y}_{1}}&
       0&
       {z}_{0} {b}_{02}+{z}_{1} {b}_{12}&
       -{z}_{0} {b}_{01}-{z}_{1} {b}_{11}\\
       {-{y}_{1}}&
       {-{y}_{2}}&
       -{z}_{0} {b}_{02}-{z}_{1} {b}_{12}&
       0&
       {z}_{0} {b}_{00}+{z}_{1} {b}_{10}\\
       {-{y}_{2}}&
       {-{y}_{3}}&
       {z}_{0} {b}_{01}+{z}_{1} {b}_{11}&
       -{z}_{0} {b}_{00}-{z}_{1} {b}_{10}&
       0\\
       \end{pmatrix}.
  $$
The projection $\pi_1$ induces a map $C^2\to \PP^1$ which is a covering of degree $3b+2$.
 \end{proposition}
 
 \begin{proof} 
 Let $I_C=\langle y_1^2-y_0y_2,y_1y_2-y_0y_3,y_2^2-y_1y_3 \rangle$ denote the homogeneous ideal
 of $C \subseteq \PP^3$. Since
 $$
 \det M \equiv -({y}_{1} {b}_{00}+{y}_{2} {b}_{01}+{y}_{3} {b}_{02}-{y}_{0}
      {b}_{10}-{y}_{1} {b}_{11}-{y}_{2} {b}_{12})^2 \mod  I_C,
 $$
the curve $C$ intersects the branch divisor of $\delta$ tangentially in $3(b+1)$ distinct points for general choices of $b_{ij}$ and the preimage of $C$ in $\PP(1,1,1,1,b+1)$ has two components defined by $I_C$ and
 $$ w \pm ({y}_{1} {b}_{00}+{y}_{2} {b}_{01}+{y}_{3} {b}_{02}-{y}_{0}
      {b}_{10}-{y}_{1} {b}_{11}-{y}_{2} {b}_{12})=0.$$
 The second statement follows by computing a primary decomposition of $I_C+\langle f \rangle \subseteq \QQ[z_0,z_1,y_0,y_1,y_2,y_3,b_{00},\ldots,b_{12}]$, where $f$ is given as in Proposition \ref{pro:Xb1}(a).
 \end{proof}

\subsection{The second small resolution}
Finally, we show that the variety $X_b^2$ defined in \S\ref{subsec:overview} is another small resolution of $Y_b$, and we finish the proof of the main theorem.

\begin{proposition}\label{pro:Xb2}
The variety $X_b^2$ is another small resolution of $Y_b$.
\end{proposition}

\begin{proof} 
As in the proof of Proposition \ref{pro:A1}, write 
$$M=\begin{pmatrix}
a_0 & a_1 \cr
a_1 & a_2 \cr
\end{pmatrix}
$$ 
for forms $a_i$ of degree $b+1$ on $\PP^3$ as in \eqref{eq:M}. Let $(u_0:u_1)$ be the coordinates on $\PP^1$. Consider the subvariety of $\PP^1\times \PP(1,1,1,1,b+1)$ defined by the $2 \times 2$ minors of the matrix
$$
\begin{pmatrix}
{a}_{0}&
{a}_{1}-w\\
{a}_{1}+w&
      {a}_{2}\\
  u_1& -u_0\\    
      \end{pmatrix};
$$
compare to \eqref{eq:firstsmall}. This is another small resolution of $Y_b$, and we will show that it is isomorphic to $X_b^2$, as defined in \S\ref{subsec:overview}. To this end, it suffices to show that the base locus of the linear system $|\sO_{X_b^1}(-1,b+1)|$ is precisely the collection of the $(b+1)^3$ exceptional curves of the small resolution $\xi_1\colon X_b^1 \to Y_b$, see Proposition \ref{pro:A1}. 

We have $\zeta^{-1}\big(\{u_1=0\}\big)=V(a_0,a_1+w,w^2+\det M)$. In $X_b^1$ this fiber is contained in $V(a_0,f)$. Since 
$$
f \equiv z_1(2z_0a_1+z_1a_2) \mod a_0
$$
is reducible, the locus $V(a_0)$ cuts $X_b^1$ in two components: $V(z_1)\in|\sO_{X_b^1}(1,0)|$ and 
$$V(a_0,2z_0a_1+z_1a_2)\in|\sO_{X_b^1}(-1,b+1)|.$$
By analysing $\zeta^{-1}\big(\{u_0=0\}\big)$, we get that another divisor in this linear system is $V(a_2,z_0a_0+2z_1a_1)$. Hence, the base locus of $|\sO_{X_b^1}(-1,b+1)|$ is the zero locus $V(a_0,a_2,2z_0a_1,2z_1a_1)=V(a_0,a_2,a_1)$, which is precisely the collection of the $(b+1)^3$ exceptional curves of $\xi_1$.
\end{proof}

Finally, our main result follows from combining all these results with the following theorem.

\begin{theorem}
The Picard group of $X_b^1$ is $\Pic(X_b^1) \simeq \Pic(\PP^1\times \PP^3)$. The nef, effective and movable cones of $X_b^1$ are
$$\Nef(X_b^1) =\langle (1,0), (0,1) \rangle $$
and 
$$\Eff(X_b^1)=\Mov(X_b^1) = \langle  (1,0), (-1,b+1) \rangle. $$
The variety $X_b^1$ has precisely two minimal models and one nontrivial birational automorphism $\iota$ of order two. The automorphism $\iota$ interchanges the two models.
\end{theorem}

\begin{proof}
The isomorphism $\Pic (X_b^1) \simeq \Pic(\PP^1\times \PP^3)$ follows from $H^2(X_b',\sO_{X_b'})=0$ and from $H^2(\PP^1 \times \PP^3,\ZZ) \simeq H^2(X_b^1,\ZZ)$, see \cite[\S3.2.A]{Laz04}.

We first prove that $\Nef(X_b^1)= \langle (1,0), (0,1) \rangle.$  Indeed, the fibres $\PP^1$ of the small resolution $\xi_1\colon X_b^1\to Y_b$ have intersection number $0$ with $\sO_{X_b^1}(0,1)$ and $1$ with $\sO_{X_b^1}(1,0)$. Thus, $\sO_{X_b^1}(\alpha, \beta)$ with $\alpha<0$ has negative intersection number with these curves. On the other hand, the curves which arise as the intersection of a fiber of $\pi_1\colon  X_b^1 \to \PP^1$ with $\pi_3^{-1}(H)$, where $H \in |\mathcal O_{\PP^3}(1)|$, have positive intersection number with $\sO_{X_b^1}(0,1)$ and intersection number $0$ with $\sO_{X_b^1}(1,0)$. Since these curves form a covering family, the line bundles $\sO_{X_b^1}(\alpha,\beta)$ with $\beta<0$ are neither nef nor effective.

Next we compute the effective and movable cone. Since $\sO_{X_b^1}(-1,b+1)$ has no fixed component by the proof of Proposition \ref{pro:Xb2}, we have $\langle (1,0), (-1,b+1) \rangle \subseteq \Mov(X_b^1)$. To see that this coincides with $\Eff(X_b^1)$ we note that the two small resolutions $X_b^1$ and $X_b^2$ of $Y_b$ coincide in codimension $1$ and are isomorphic as abstract varieties. Thus, we have
$$
h^0\big(X_b^1,\sO_{X_b^1}(\alpha,\beta)\big)=h^0\big(X_b^2,\sO_{X_b^2}(\alpha,\beta)\big) =h^0 \big(X_b^1,\sO_{X_b^1}(-\alpha,\alpha(b+1)+\beta)\big).
$$
In particular, these groups are zero for $\alpha>0$ and $\beta<0$ and
$$
\Eff(X_b^1) =\Mov(X_b^1) =  \langle (1,0), (0,1) \rangle \cup \langle  (0,1), (-1,b+1)\rangle. 
$$
The interiors of the two subcones are ample on $X_b^1$ and $X_b^2$, respectively.
\end{proof}

All computations in Macaulay2 can be found in \cite{LS19}.

\subsection*{Supporting file}
A supporting file for this paper is available on the journal website.

	\bibliographystyle{amsalpha}
	\bibliography{biblio}

\end{document}